\documentclass[titlepage,draft,12pt]{article}
\usepackage{theorem}
\usepackage{array}
\usepackage[reqno]{amsmath}
\usepackage{enumerate}
\usepackage{fancybox}
\usepackage{color}
\input epsf.sty
\usepackage[final]{graphics}

\usepackage{url}

\theoremstyle{change}
\newtheorem{proclaim}{PROCLAIM}[section]
\newtheorem{theorem}[proclaim]{Theorem}
\newtheorem{definition}[proclaim]{Definition}
\newtheorem{lemma}[proclaim]{Lemma}

\newtheorem{corollary}[proclaim]{Corollary}
\newtheorem{example}[proclaim]{Example}

\outer\def\proclaim #1. #2\par{\medbreak \noindent{\bf#1.\enspace}{\sl#2}\par
  \ifdim\lastskip<\medskipamount
  \removelastskip\penalty55\medskip\fi}
\def\state #1. { \noindent{\bf#1.\enspace}}

\newcommand{\rwtxt}[1]{\,\;\hbox{#1}\;\,}
\newcommand{\Text}[1]{\,\quad\hbox{#1}\quad\,}

\newcommand{\argmax}{{\rm argmax}}
\newcommand{\argmaxinf}{{\rm argmaxinf}}
\newcommand{\argmin}{{\rm argmin}}

\newcommand{\dom}{{\rm dom}}
\newcommand{\comp}{\,{\raise 1pt \hbox{$\scriptstyle\circ$}}\,}
\newcommand{\eqhyph}{\hbox{\rm -}}

\newcommand{\reals}{{I\kern-.35em R}}
\newcommand{\Reals}{\overline{I\kern-.35em R}}
\newcommand{\natnums}{{{\rm l} \kern -.13em {\rm N} }}
\newcommand{\nats}{{I\kern -.35em N}}
\newcommand{\snats}{{I\kern -.29em N}}
\newcommand{\rats}{{Q\kern -.64em \raise 1pt \hbox{$\scriptstyle |$}\;\,}}
\newcommand{\srats}
	{{Q\kern -.56em \raise 1.2pt \hbox{$\scriptscriptstyle /$}\,}}
\newcommand{\ints}{Z\kern -.46em Z}
\newcommand{\ball}{{I\kern -.35em B}}
\newcommand{\pluss}{\hskip1pt \raise1pt\vbox{\hrule width6pt \vskip1pt \hrule
                    width6pt} \kern-4pt{\lower1pt\hbox{\vrule height6pt
		    \kern1pt\vrule height6pt}}\hskip5pt}
\newcommand{\eop}
	{\hfill{$\vcenter{\hrule height1pt \hbox{\vrule width1pt height5pt 
   	 \kern5pt \vrule width1pt} \hrule height1pt}$} \medskip}

\newcommand{\linfty}{{\scriptscriptstyle \infty}}

\newcommand{\lplus}{{\scriptscriptstyle +}}

\newcommand{\setd}{{ d \kern -.15em l}} 
\newcommand{\hatsetd}{ d \hat{\kern -.15em l }}

\renewcommand{\epsilon}{\varepsilon}
\renewcommand{\phi}{\varphi}

\newcommand{\lset}{\big\lbrace}
\newcommand{\mset}{{\,\big\vert\,}}
\newcommand{\rset}{\big\rbrace}

\newcommand{\Mset}{{\,\Big\vert\,}}

\newcommand{\tto}{\;{\lower 1pt \hbox{$\rightarrow$}}\kern -12pt
           \hbox{\raise 2.5pt \hbox{$\rightarrow$}}\;}
\newcommand{\overto}[1]{\,{\raise 0pt\hbox{$\rightarrow$}}\kern -9pt
     \hbox{\lower 3pt \hbox{$\scriptscriptstyle#1$}}\hskip6pt}
\newcommand{\underto}[1]{\,{\lower 1pt\hbox{$\rightarrow$}}\kern -9pt
     \hbox{\raise 4pt \hbox{$\,\scriptscriptstyle#1$}}\hskip7pt}
\newcommand{\bigoverto}[1]{{\raise 0pt\hbox{$\,\longrightarrow$}}\kern -16pt
     \hbox{\lower 3pt \hbox{$\scriptscriptstyle#1$}}\hskip4pt}
\newcommand{\bigunderto}[1]{\,{\lower 1pt\hbox{$\longrightarrow$}}\kern -16pt
     \hbox{\raise 4pt \hbox{$\,\scriptscriptstyle#1$}}\hskip6pt}
\newcommand{\bigbigto}[2]{\,{\raise 0pt\hbox{$\,\longrightarrow$}}\kern -16pt
     \hbox{\lower 3pt \hbox{$\scriptscriptstyle#2$}}\kern -10pt
     \hbox{\raise 4pt \hbox{$\,\scriptscriptstyle#1$}}\hskip7pt}
\newcommand{\downto}{{\raise 1pt \hbox{$\scriptstyle \,\searrow\,$}}}
\newcommand{\upto}{{\raise 1pt \hbox{$\scriptstyle \,\nearrow\,$}}}

\newcommand{\notimply}
	{\quad\hbox{$\Longrightarrow \kern -14pt {/}$}\hskip6pt\quad}
\newcommand{\rightrightarrows}{\,{{\lower4pt\hbox{$\scriptstyle\rightarrow$}}
    \atop{\raise8pt\hbox{$\scriptstyle\rightarrow$}}}\,} %supplanted by \tto

\def\eto{\,{\lower 1pt\hbox{$\rightarrow$}}\kern -11pt
     \hbox{\raise 4pt \hbox{$\, \scriptstyle e$}}\hskip7pt}
\def\hto{\,{\lower 1pt\hbox{$\rightarrow$}}\kern -11pt
     \hbox{\raise 4pt \hbox{$\, \scriptstyle h$}}\hskip7pt}
\def\gto{\,{\lower 1pt\hbox{$\rightarrow$}}\kern -11pt
     \hbox{\raise 4.5pt \hbox{$\, \scriptstyle g$}}\hskip7pt}
\def\pto{\,{\lower 1pt\hbox{$\rightarrow$}}\kern -11pt
     \hbox{\raise 4.5pt \hbox{$\, \scriptstyle p$}}\hskip7pt}
\def\cto{\,{\lower 1pt\hbox{$\rightarrow$}}\kern -11pt
      \hbox{\raise 4pt \hbox{$\, \scriptstyle c$}}\hskip7pt}
\newcommand{\awto}{\,{\lower 1pt\hbox{$\longrightarrow$}}\kern -16pt
     \hbox{\raise 4pt \hbox{$\,\scriptstyle aw$}}\hskip6pt}

\newcommand{\low}[1]{{\lower1pt \hbox{$\scriptstyle #1$}}}
\newcommand{\loww}[1]{{\lower2pt \hbox{$\scriptstyle #1$}}}
\newcommand{\high}[1]{{\raise1pt \hbox{$\scriptstyle #1$}}}

\newcommand{\cI}{{\cal I}}

\newcommand{\nsum}{\mathop{\sum}\nolimits}
\newcommand{\nprod}{\mathop{\prod}\nolimits}

\newcommand{\ninf}{\mathop{\rm inf}\nolimits}
\newcommand{\nsup}{\mathop{\rm sup}\nolimits}
\newcommand{\nmax}{\mathop{\rm max}\nolimits}

\newcommand{\nargmax}{\mathop{\rm argmax}\nolimits}

\newcommand{\bfxi}{\mbox{\boldmath $\xi$}}

\newcommand{\sothat}{\mathop{\rm so\;\, that\;}\nolimits}
\newcommand{\suchthat}{\mathop{\rm \;such\; that\;}\nolimits}

\newcommand{\lwdy}[2]{\mathrel{\mathop
        {\raisebox{0.1ex}{\null$#1$}}{\hbox{\kern -1.0em
	{\raisebox{-0.8ex}{$\scriptstyle{\;\to #2}$}}}}}}
\newcommand{\lwwdy}[2]{\mathrel{\mathop
        {\raisebox{0.2ex}{\null$#1$}}{\hbox{\kern -1.0em
	{\raisebox{-1.1ex}{$\scriptstyle{\;\to #2}$}}}}}}
\newcommand{\slwwdy}[2]{\scriptsize{{\mathrel{\mathop
        {\raisebox{0.2ex}{\null$#1$}}{\hbox{\kern -1.0em
	{\raisebox{-1.1ex}{$\scriptstyle{\;\to #2}$}}}}}}}}
\newcommand{\cc}{\!:}

\newcommand{\Ex}{{I\kern-.35em E}}

\begin{document}
\bibliographystyle{plain}
\begin{titlepage}
\vglue 0.5cm
\begin{center}
\begin{large}
{\bf SOLVING DETERMINISTIC AND STOCHASTIC\\
EQUILIBRIUM PROBLEMS VIA AUGMENTED WALRASIAN 
%\footnote{Julio Deride's work was initially supported by CMM =
%Centro de Modelamiento Matem{\'a}tico, Universidad de Chile.  Later,
%this material is based upon work supported in part by the
%U. S. Army Research Laboratory and the U. S. Army Research Office
%under grant numbers 00101-80683, W911NF-10-1-0246 and W911NF-12-1-0273.}
}

\smallskip
\end{large}
\vglue 3.0truecm
\begin{tabular}{lclcl}
  \begin{large} {\sl Julio Deride %% AUTHOR #1
                                  } \end{large} & \ \ {\phantom{\&}} \ \ &
  \begin{large} {\sl Alejandro Jofr{\'e} %% AUTHOR #2
                                  } \end{large} & \ \ {\phantom{\&}} \ \ &
  \begin{large} {\sl Roger J-B Wets %% AUTHOR #3
                                  } \end{large} \\
  \\
  Mathematics  && CMM \& DIM &&   Mathematics \\
  Univ. California, Davis && Universidad de Chile && Univ. California, Davis \\
  jderide@math.ucdavis.edu && ajofre@dim.uchile.cl && rjbwets@ucdavis.edu 
\end{tabular}
\end{center}
\vskip 1.5truecm
\noindent {\bf Abstract}. \quad We described a method to solve
deterministic and stochastic Walras equilibrium models based on
associating with the given problem a bifunction whose maxinf-points
turn out to be equilibrium points. The numerical procedure relies on an
augmentation of this bifunction. Convergence of the proposed 
procedure is proved by relying on the relevant lopsided 
convergence. In the dynamic versions of our models, deterministic 
and stochastic, we are mostly concerned with models that equip the 
agents with a mechanism to transfer goods from one time period 
to the next, possibly simply savings, but 
also allows for the transformation of goods via production.

\vskip 1.5truecm
\halign{&\vtop{\parindent=0pt
   \hangindent2.5em\strut#\strut}\cr
\hglue 1.50cm  augmented Walrasian, progressive hedging algorithm.\hfill\break\cr
%\hfill\break
%\hglue 1.50cm  exponential convergence.  \cr\cr
{\bf Date}:\quad \ \today \cr}
%% {\bf Date}:\quad \ September 10, 2005, \quad {\bf Revised}: \quad March 12, 2006 \cr}
\end{titlepage}
\baselineskip=15pt

% \usepackage{lscape}

% April 8, 2015 (Wednesday, 11:43 PDT) rjbw

\section{Introduction}

The economic equilibrium model proposed by Arrow and Debreu \cite{ArrD54} for a
competitive economy implicitly assumes that the entire economic activity will take
place in a single time span, implicitly instantly.  As soon as one includes the agent's
concerns about the future, one has to consider a dynamic component and take into
account the uncertainty about this future. Of course, that's bound to enrich the model,
raising in the process a wide variety of modeling issues.  In this article we are
only going to be concerned with numerical procedures to solve trimmed down stochastic
Walras equilibrium models where goods are transferred from time 0 to time 1 via
``home production'' which includes the possibility of simple retention; the extension
of this approach to include financial markets is presently under development.

The overall approach to the stochastic case is based on a result that allows us to
proceed with the calculation of the equilibrium for each particular state, a
decomposition type approach. This means that our first task is going to be the
development of a method that will arrive at an equilibrium rather efficiently in a
deterministic, but dynamic, environment. We start with the pure exchange model of
Arrow-Debreu, next consider a two time-periods dynamic version and then proceed to
deal with the stochastic version of this model.  We rely on an augmentation method
applied to what we call the Walrasian, essentially, the function `supposedly
solved' by the Walrasian auctioneer and a maxmin characterization of an equilibrium point. 
The fact that the theory allows us to proceed, in the iterative process, with approximate 
equilibria turns out to be critical in the development of the overall numerical scheme.

Our approach deviates, even in the deterministic case, from the path-breaking methods
suggested by Scarf and Hansen \cite{ScrH73}, Eaves \cite{Eves11, BrDeEa96:eaves}, 
Saigal \cite{Sgal83}, and other approximation strategies described in the books by 
Judd \cite{Judd98}, and Brown and Kubler \cite{BrKu08}. These earlier methods are
not efficient when the economies have a large number of goods or agents, even for reaching approximate equilibria.
Moreover, in stochastic environments these results are prohibitively
time-expensive.  

In this paper we develop an approach based on an augmented Walrasian technique and a lopsided convergence approximation procedure, which allows us to cope with large equilibrium problems including uncertainty and heterogeneity on the agents. By using this approach we have designed a two-phase algorithm without computing derivatives of the demand function. We report several numerical experiments for equilibrium problems involving up to 5 agents, 7 goods and 10 stochastic scenarios, which can be easily expanded in number of agents and goods. Finally, the procedure proposed in this paper might be parallelized in terms of agents and the multi-start strategy.

% Ago 04, 2014 (Monday, 17:52 PST) jd

\section{The Arrow-Debreu model} \label{arrowDebreu}

To set the stage and fix terminology and notation, let's start with
the barter, or pure exchange model, of Arrow-Debreu \cite{Dbru59}.
A finite number of (individual) agents $i\in\cI$ with initial endowments
$\lset e_i \in \reals^L, \, i\in \cI \rset$,
consisting of a finite number of goods, to be bartered so as to
maximize, individually, their upper semicontinuous (usc) concave
utility functions $\lset u_i: \reals^L \to
[-\infty, \infty), \, i\in\cI \rset$ that depend on the level of the
acquisitions
$x_i(p)\in \reals^L$ of these goods, potentially for ``consumption'';
one refers to $X_i = \dom \, u_i = \lset x\in \reals^L
\mset u_i(x) > -\infty \rset$ as the {\em survival sets}; note that
the concavity of $u_i$ implies that the survival set
$X_i$ is convex, typically unbounded. The value
to assign to each good, in this trading process, depends on a
{\em market price system}
$0 \neq p\in \reals_\lplus^L$ that will restrict each agent to limit
the ``market value'' of its acquisitions to the ``market value'' of its
endowment, i.e., $\langle p,x \rangle \leq
\langle p, e_i \rangle$; since these prices don't necessarily reflect
monetary prices, the ``values'' are often referred to as
{\em units of account}.  Given $p\in \reals_\lplus^L$,
each agent maximizes its utility subject to its budgetary constraint,
i.e.,
\[ x_i(p) \in \nargmax_{x\in X_i\subset\reals_\lplus^L} u_i(x) \]
For the market to be in {\em equilibrium} the total demand must not
exceed total supply, i.e., with $s$ designating the {\em excess supply
function},
\[ s(p) = \nsum_{i\in\cI} (e_i - x_i(p)) \geq 0. \]
Since, we haven't ruled out the possibility that at equilibrium the
prices of some goods might turn out to be 0, one can also write this
condition in terms of a geometric variational inequality:
\[  - \nsum_{i\in\cI} (e_i - x_i(p)) = -s(p) \in N_{\reals_\lplus^L}(p), \]
where $N_C(p)$ denotes the {\em normal cone} of variational analyis
\cite{VaAn} to the set $C$ at $p$, 
or still, must be solutions of the {\em linear complementarity problem},
\[ 0 \leq p \perp \nsum_{i\in\cI} (e_i -d_i(p)) \geq 0. \]

\noindent 
Since the {\em budgetary constraints} are positively homogeneous in $p$
and $p \neq 0$, no additional restriction is introduced by insisting
that the price system should be scaled so that it lies in the unit
simplex $\Delta^L = \lset p\in \reals_\lplus^L \mset \langle p, e \rangle
= 1 \rset$.  This is often included in the formulation of the problem to
enable appealing to a fixed point argument to establish existence or to
provide boundedness in the design of a computational scheme.\\
%\footnote{Depending on the aims being pursued one might resort to other
%scaling strategies, for example by insisting that the price of a specific
%(normalizing) good be 1.  When fiat money is included in our goods' bundle,
%it can be used as a ``num{\'e}raire'' and then the price system can be
%normalized so that all other prices are then expressed in terms of the
%fiat money denomination.}.\\

%%%%%%%%%%%%%%%%%%%%%%%%%%%%%%%%%%%%%%%%%%%%%%%%%%%%%%%% jd
\noindent
Additionally, we introduce a natural bound for each agent demand 
function $x_i(p)\leq \sum_{i\in\cI} e_i$  (see \cite[Ch.5]{Dbru59}), i.e., no agent can demand more
quantity of each good than the total amount available in the economy,
which in turn obtains a bound for $s(p)>-\infty$.

%%%%%%%%%%%%%%%%%%%%%%%%%%%%%%%%%%%%%%%%%%%%%%%%%%%%%%%%

% April 13, 2015 (Monday, 12:13 PDT) jd

\section{Augmented Walrasian}

Our assumptions, and notation introduced in the previous section,
follow those of the article ``Continuity
properties of Walras equilibrium points'' \cite{JfrW02:wlrs} which introduced the 
{\em Walrasian} function associated with this problem
\[ W(p,q) = \langle q,s(p) \rangle \;\rwtxt{on} \Delta\times\Delta \]
where $p,q \in$ the (unit) price simplex $\Delta\subset \reals^n$, and
$s$ is our excess supply function as defined in the previous section. 
Moreover, the following lemma provides that
every max-inf point of $W$ is an equilibrium price, i.e., 
if $\inf_{q\in \Delta} W(\bar{p},q)=\sup_{p\in\Delta} \inf_{q\in\Delta} W(p,q)$ 
then $s(\bar p)\geq 0$  \cite{JfrW02:wlrs}.

\begin{lemma}(\rm {Walras equilibrium prices and maxinf-points}).
\label{lemm:maxinf}
Every maxinf-point $\bar{p}\in\Delta$ of the Walrasian function $W$ such that 
$W(\bar{p},\cdot) \geq 0$ on $\Delta$ is an equilibrium point.
\end{lemma}

\state Proof.
If $\bar{p}$ is a maxinf-point of the Walrasian with $W(\bar{p},\cdot) \geq 0$, 
it follows that for all unit vectors $e^j=(0,\ldots,1,\ldots)$, the $j$-th 
entry is 1, $\langle e^j, s(\bar{p}) \rangle \geq 0$ which implies $s(\bar{p}) \geq 0$.
\eop

\noindent
The condition of $W(\bar{p},\cdot) \geq 0$ follows by the definition 
of $W$ and noting that for every price $p\in\Delta$, and under local
nonsatiation preferences assumption, $W(p,p)=0$. Furthermore, the 
converse of lemma (\ref{lemm:maxinf}) also holds, i.e., every equilibrium 
point is a maxinf-point of $W$ \cite[Prop.2.4]{JfrW11:var_conv_mot}.\\

\noindent
Existence of equilibrium prices can be seen as the existence of
max-inf points for the corresoponding Walrasian. Under general 
conditions for the upper-semicontinuity of the excess supply function,
it is easy to see that $W$ is a Ky Fan function \cite{Aub06:app_nla}, and the 
existence of max-inf point is provided by \cite[Theorem 6.3.5]{Aub06:app_nla}\\

\noindent
Since our basic approach, first suggested by A. Bagh
\cite{Bagh02}, is related to that for the augmented Lagrangian,
it's informative to consider the bifunction that might
have led to the Walrasian in a standard non-convex duality scheme
\cite[\S 11.K]{VaAn}.  Let's introduce a {\em pre-Walrasian}
obtained as a restricted-partial conjugate, with respect to the
$q$-variable, i.e., for all $p\in \Delta$ 
\[V(p,u) = \nsup_{\;z\in\Delta}\,[\,\langle u,z\rangle - W(p,z)\,].\]
$V(p,\cdot)$ is clearly convex and one can think of the family
of bifunctions $\lset V(\cdot,u),\,u\in \reals^n \rset$ as 
`perturbations' of a `{\em fundamental}' primal-problem 
\[ \rwtxt{find}\;\hat p \in\nargmax_{\;p\in\Delta}\,v(p) \Text{where}
  v(p) = V(p,0) = -\inf_{\,q\in\Delta}\,[\,W(p,q)\,]. \]
By conjugacy, since the functions $q\mapsto W(p,q)$ on $\Delta$
are proper, lower semicontinuous (lsc) and convex, so 
are the functions $u\mapsto V(p,u)$.
%; $\iota_C$ is the indicator function of the set $C$.
Note that $\min_{\,q\in\Delta}\,\langle q,s(p) \rangle$ will yield the
$q$ that generates the smallest convex combination of the elements of
$s(p)$. So, if for any $l$, $s_l(p) <0$, it follows that $v(p)
>0$. Thus, $\hat p$ will be such that an element of the vector $s(\hat
 p)$ will be as negative as possible it will minimize the
 $\ell^\linfty$-norm of $s(p)$. \\  

\noindent
The process of going from $v$ to the collection $\{V(\cdot,u),
u\in \reals^n\}$ is well-understood; it can be viewed as associating to
a particular optimization problem, $\max \lset v(p) \mset p\in\Delta
\rset$, a perturbed collections that leads to the analysis of 
stability. But in our setting what is this particular
optimization problem? It can be viewed as the {\em Walrasian auctioneer's
problem}. It's easy to see that it's optimal value is 0 which is
attained when the Walrasian auctioneer has selected a price system that
yields an equilibrium. Generally it's not a concave function, and
certainly not a strictly concave function, and thus one can't expect a
unique maximizer which, precisely, correspond to the well-know fact
that, in general, Walras equilibrium points are not unique. \\

\noindent
In order to compute equilibrium points for an economy, we propose
a strategy to find a max-inf point of $W$ by an approximating scheme. 
Our first goal is to build a family of approximating bifunctions 
by relying on an augmentation technique.
Let $\sigma\cc \reals^n\to \Reals$ be an augmenting function, i.e.,
it's convex, $\argmin \sigma = \{0\}$ and $\min \sigma = 0$. Typically,
$\sigma = |\,\cdot\,|$ is chosen to be a norm but depending on the
application it could be quite different; recall that we can even choose
$\sigma$ to take on the value $\infty$, for example, it could be a norm
of some type restricted to a ball centered at 0, or even more exotic. \\

\noindent 
Given the augmenting function $\sigma$ and a scalar $r>0$, the 
{\em augmented Walrasian}, by definition, is
\[\Tilde W_r(p,q) = \nsup_{\;u\in \reals^n}\, \lset \langle q,u \rangle
	- V(p,u) - r\sigma(u) \rset. \]
%%%%% We shall rely extensively on the following identity:
\noindent
For a fixed $p$, considering the convexity of $V(p,\cdot)$
and $\sigma$, one can re-write the definition of $\Tilde W_r$
as a partial conjugate w.r.t the $u$-variable. Additionally, by 
the property of conjugation of a sum and the definition
of the epi-sum ($\pluss$), we have the following chain of identities
\begin{eqnarray*}
\Tilde W_r(p,q) & = & \left(V(p,\cdot)+r\sigma\right)^*(q)\\
& = & {\rm cl}\,\lset \left((V(p,\cdot)^*\pluss(r\sigma)^*\right)(q)\rset\\
& = & \ninf_{\;z}\, \lset (V(p,\cdot))^*(q-z) + (r\sigma)^*(q) \rset\\
& = & \ninf_{\;z}\, \lset W(p,q-z) + r\sigma^*(r^{-1}z) \rset
\end{eqnarray*}
where $\sigma^*$ is the conjugate of $\sigma$, i.e., 
$\sigma^*(v)=\sup_x\lset \langle v,x\rangle -\sigma(x) \rset$. Thus,
the augmented Walrasian function in its final form is the infimum 
of a convex function, and depending of our choice of 
$\sigma$, possibly quadratic or linear.\\

%%%%%%%%%%%%%%%%%%%%%%%%%%%%%%%%%%%%%%%%%%%%%%_beg jd 
%see the section on Augmented Lagrangians in Rockafellar \& Wets,
%{\it Variational Analysis}, Exercise 11.56. \\

%\noindent
%From now on, it will be assumed that $\sigma = |\,\cdot\,|$, is a
%norm whose dual norm will be denoted by $|\,\cdot\,|_o$. Then, 
%one can express the augmented Walrasian in the following form:
%\[\Tilde W_r(p,q) = \ninf_{\;z\in \reals^n}\, \Lset W(p,z) \Mset
%	| z-q|_o \leq r. \Rset \]
%This is the form, at least at present, that seems to be potentially the best for
%computational purposes, but that's not been demonstrated since no other
%choice of augmenting function $\sigma$ has been tried. 
%The natural candidates for $|\,\cdot\,|$ are
%the $\ell^2$- and $\ell^1$-norms that lead respectively to
%$|\,\cdot\,|_o$ as the indicator functions of the (euclidean) unit ball
%and the $\ell^\linfty$-ball. The ball associated with the $|\,\cdot\,|_o$
%will be denoted by $\ball_o$. \\

\noindent
Additionally, we establish a definition for an approximating 
equilibrium point, given by a price such that the associated excess
supply function is close to satisfy the market clearing condition. 
More precisely, the definition can be stated as follows: 

\noindent
\begin{definition}{\rm (approximate maxinf-points).}
For $\epsilon \geq 0$, $p_\epsilon$ is said to be an 
\emph{approximate equilibrium point} or \emph{approximate maxinf-point} 
of $W$ if $ \inf W(p_\varepsilon,\cdot) \geq \sup\inf W-\varepsilon$, 
and the set of all such approximating maxinf-points is denoted by 
$\varepsilon \eqhyph \argmaxinf\, W$.
\end{definition}

\noindent
Note that given an approximating equilibrium price, 
$p_\epsilon$, one can adjust the agents's intial endowments by a 
fraction of $\epsilon$ and make $p_\epsilon$ and equilibrium
price.\\

\noindent
The next step is to establish a connection between the convergence of 
augmented Walrasian approximating equilibrium points and the goal of finding 
an equilibrium price for the original economy. Considering the 
family of augmented Walrasian perturbations, and a sequence of their 
corresponding approximate max-inf points, one should be able to guarantee 
a convergence result of this sequence of points, given the convergence of the 
family of augmented functions. This condition can be obtained by appealing to
\emph{lopsided convergence}, or \emph{lop-convergence}, of the augmented 
Walrasian to the Walrasian. Given the compactness 
of the domain, one doesn't have to appeal to the (more comprehensive) definition
of lopsided convergence it suffices to refer to a more restrictive version,
namely \emph{tight} lopsided convergence; for the general definition and
further details, consult \cite{JfrW09:var_conv,JfrW11:var_conv_mot}.

\begin{definition}{\rm (tight lopsided convergence).}
A sequence in finite-valued bivariate functions, fv-biv($\reals^{n+m}$), defined 
over a compact set $C\times D$, $\lset F^\nu:C \times D \to \reals \rset_{\nu \in \nats}$ 
{\em lop-converges tightly} to a function $F:C\times D\to \reals$, 
also in fv-biv($\reals^{n+m}$), if 

\begin{itemize}
\item[(a)] for all $y\in D$, and all $(x^\nu\in C)\to x\in C$, 
there exists $(y^\nu \in D) \to y$ such that 
$$\limsup_\nu F^\nu(x^\nu,y^\nu) \leq F(x,y),$$

\item[(b)] for all $x\in C$, there exists $(x^\nu \in C) \to x$ 
such that given any $(y^\nu \in D) \to y\in D$, 
$$\liminf_\nu F^\nu(x^\nu,y^\nu)\geq F(x,y),$$
\end{itemize}

\end{definition}

\noindent The desired convergence result for the equilibrium points 
follows from adaptating the tight lop-convergence given by 
\cite[Theorem 3.2]{JfrW11:var_conv_mot}, to our case with compact 
(and invariant) domains.

\begin{theorem}
\label{mot_3.2}{\rm (convergence of maxinf-points, \cite[Theorem 3.2]{JfrW11:var_conv_mot})}.
Let $C\times D$ be a compact subset of $\reals^{n+m}$. When the bifunctions 
$\lset F^\nu\rset_{\nu\in\nats}$ lop-converge tightly to $F$, 
all in fv-biv($C\times D$) with $\sup\inf F$ finite, and 
$\epsilon^\nu \downto \epsilon \geq 0$, then every cluster point $\bar{x}\in C$
 of a sequence of $\epsilon^\nu$-maxinf points of the bifunctions $F^\nu$ is 
an $\epsilon$-maxinf point of the limit function $F$.

In particular, this implies that in these circumstances, every cluster point 
of a sequence of maxinf-points of the bifunctions $F^\nu$ is a maxinf-point 
of the lop-limit function $F$.
\end{theorem}

\noindent
In order to obtain our convergence result for approximating maxinf points, the 
following result is an application of the previous theorem in our framework. It
tell us that tight lopsided convergence of the augmented Walrasian entails 
convergence of equilibrium points.

\begin{theorem}{\rm (convergence of $\epsilon$-maxinf points).}\label{conv_mi_ad}
Suppose that $p \mapsto s(p)$ is usc on $\Delta$. Consider the non-negative sequences 
$\lset r^\nu,\, \nu \in \nats \rset$ and $\lset \epsilon^\nu,\, \nu \in \nats \rset$ 
such that $r^\nu \upto \infty$, $\epsilon^\nu \downto \epsilon\geq 0$. 
Let $\lset W^{\nu},\, \nu \in \nats \rset$ be a family of augmented 
Walrasian functions associated which each augmenting parameter $r^\nu$. Let 
$p^{\nu} \in \epsilon^{\nu} \eqhyph \argmaxinf\, W^{\nu}$ 
and $\bar{p}$ be any cluster point of $\lset p^{\nu},\,nu \in \nats\rset$.
Then $\bar{p} \in \epsilon \eqhyph \argmaxinf\, W$.
\end{theorem}

\state Proof.
It suffices to show that $\lset W^\nu,\, \nu \in \nats \rset$ lop-converges 
tightly to $W$ and conclude by Theorem (\ref{mot_3.2}) the convergence 
of a (sub)sequence of $\epsilon^\nu$-maxinf points. In order to prove 
tight lopsided convergence, let $q\in \Delta$, $\lset p^\nu,\, \nu \in \nats \rset \to p\in \Delta$. 
Define $q^\nu \equiv q,\,\nu \in \nats$. Then
\[ W^\nu (p^\nu, q^\nu)=\ninf_{\;z\in \reals^n}\, \lset W(p^\nu,z) +
	r^\nu*\sigma^*(q^\nu-z) \rset \leq W(p^\nu,q^\nu), \]
and as the function $p \mapsto s(p)$ is usc, 
\[ \limsup W^\nu (p^\nu, q^\nu) \leq \limsup W(p^\nu,q) \leq W(p,q). \]

\noindent
On the other hand, let $p\in \Delta$ and $\lset q^\nu,\, \nu \in \nats \rset \to q$. 
By compactness of $\Delta$, $q\in \Delta$ and defining 
$p^\nu = p,\, \nu \in \nats$, $W^\nu(p,q)$ is the inf-projection 
of the function $F^\nu(q,z) = W(p,q-z) + r^\nu*\sigma^*(z)$ 
in the $z$-variable. Thus, $F^\nu$ is level bounded in $z$ 
locally uniform in $q$ and therefore $W^\nu(p,\cdot)$ is lsc 
by \cite[Theorem 1.17]{VaAn}. Finally, 
       
\[ \liminf W^\nu(p^\nu,q^\nu) \geq W(p,q), \]  

\noindent since for any $q_0 \in \Delta$, 
$W^{\nu}(p,q_0)\to W(p,q_0)$ as $\nu \to \infty$ and the conclusion
follows from a standard diagonal argument.
\eop

\noindent The following inmediate corollary of this theorem, 
with $\epsilon = 0$, plays a pivotal role form a numerical
viewpoint

\begin{corollary}{\rm ($\epsilon$-maxinf and equilibrium points).}\label{cor_eq_conv}
Let $\epsilon^\nu \downto 0$. Then, every cluster point of a sequence of 
$\epsilon^\nu$-approximating equilibrium points of a sequence of augmented
 Walrasian functions is an equilibrium point for the original economy.
\end{corollary}

%%%%%%%%%%%%%%%%%%%%%%%%%%%%%%%%%%%%%%%%%%%%%%%%%%%%%%%%%_end jd

\noindent
The major thrust of the eventual algorithmic procedures is to replace
finding local near  maxinf-points of $W$ by finding a local saddle point of 
a $\Tilde W^\nu$ for $\nu$ large enough (but not too large to avoid numerical
instabilities). 
%The basic justification for proceeding in this fashion
%can be found in {\it Variational Analysis}, Theorem 11.59, identity
%11(30). 
Under this scheme, there are several options for choosing the augmenting 
function $\sigma$. For example, one can consider $\sigma=|\,\cdot\,|$ 
a norm whose dual norm will be denoted by $|\,\cdot\,|_o$, then one can express 
the augmented Walrasian as 
\[\Tilde W^\nu(p,q) = \min_z \left[\, W(p,z) \Mset z\in\ball_0(q,r_\nu)\cap\Delta\,\right], \]

\noindent
where $\ball_0(q,r_\nu)$ is the dual ball with center in $q$ and radius $r_\nu$.
Alternatively, for $\sigma$ be the self-dual function, i.e., 
$\sigma=\frac{1}{2}|\,\cdot\,|_2^2$, the augmented Walrasian takes the form
\[\Tilde W^\nu(p,q) = \min_z \left[\, W(p,z)+\frac{1}{2r^\nu}|z-q|_2^2 \,\Mset z\in\Delta\,\right].\]

\noindent
There is quite a variety of procedures for finding these near local
saddle-points. One possible procedure to solve the problem at hand
is described next:

\begin{itemize} 
\item At iteration $\nu+1$, given $(p^\nu,q^\nu)$ with $r = r_{\nu+1}$ ($\geq r_
\nu$), the Phase I (or primal) consists in solving
%\[q^{\nu+1} \in \argmin_{q\in\Delta} \,\left[\, \min_z W(p^\nu,z) \Mset
%	z\in \ball_o(q,r_{\nu+1})\cap\Delta \,\right]; \]
\[q^{\nu+1} \in \argmin_{q\in\Delta}\,\Tilde W^{\nu+1}(p^\nu,q) \]
note that the `internal' minimization is either that of a linear form on a
ball, this seems to favor $|\,\cdot\,|_o$ as the $\ell^\linfty$-norm,
or the self-dual augmenting function which yields an immediate solution.

\item How to carry out the next step will depend on the `shape' and
the properties of the demand functions. For example, this turns out to
be rather simple when the utility functions are of the Cobb-Douglas
type, defining the Phase II (or dual) as finding
%\[p^{\nu+1} \in \argmax_{p\in\Delta} \,\left[\, \min_z W(p,z) \Mset
%        z\in \ball_o(q^{\nu+1},r_{\nu+1})\cap\Delta \,\right]. \]
\[p^{\nu+1} \in \argmax_{p\in\Delta} \,\Tilde W^{\nu+1}(p,q^{\nu+1})\,\]
\end{itemize} 

\noindent
In virtue of the corollary (\ref{cor_eq_conv}), we know that as $r_\nu \upto \infty$,
$p^\nu \to \bar p$ a maxinf-point of $W$, equivalently an equilibrium
price system for Walras' problem. The strategy for increasing $r_\nu$
should take into account (i) numerical stability, i.e., keep $r_\nu$ as
small as possible and (ii) efficiency, i.e., increase $r_\nu$ sufficiently
fast to guarantee accelerated convergence.

%%%%%%%%%%%%%%%%%%%%%%%%%%%%%%%%%%%%%%%%%%%%%%%%%%%%_be jd

\subsection{Numerical implementation for the Arrow-Debreu model.}

\noindent
The proposed algorithm was implemented in Pyomo (Python Optimization Modeling 
Objects, \cite{HLWD12:pyomo}), a mathematical programming language based on 
Python. The problems that we solve come with the following features:
\begin{itemize}
\item In order to describe the economy, we consider utility functions 
of Cobb-Douglas and Constant Elasticity of Substitution (CES) type, and 
strictly positive aggregated initial endowment for every good.
\item For the selection of the augmenting function $\sigma$, we primarily considered 
the self-dual type, given by $\sigma=\frac{1}{2}|\,\cdot\,|_2^2$.  
\item For the agent problem, everyone has to maximize a concave utility function over 
a linear constrained set determined by budgetary constraint and nonegativity 
of the solution. This problem is solved using the interior point method, Ipopt, 
implemented by \cite{WaBi06:ipopt} (which gives satisfactory results for 
problems of this nature).
\item Phase I consists of the minimization of a quadratic objective 
function over the simplex of prices. This is solved using Gurobi solver 
\cite{gurobi}, a state-of-the-art and efficient algorithm.
\item Phase II is the critical step of the entire augmented Walrasian 
algorithmic framework. We need to overcome the (typical) lack of concavity of the 
objective function. Thus, the maximization is done without considering 
first order information and relying on BOBYQA algorithm \cite{Pow09:bobyqa}. 
which performs a sequentially local quadratic fit of 
the objective functions, over box constraints, and solves it using a 
trust-region method.
\noindent All the examples were run on a 3.30 GHz Intel Core i3-3220 processor with 4
GB of RAM memory, under Ubuntu 12.04 operating system.
%Numerical examples were run in a desktop machine, Intel® Core™ i3-3220 CPU @ 3.30GHz × 4 
%4 gb ram, under Ubuntu 12.04
\end{itemize}

In what follows, a set of numerical examples is described.
The first example corresponds to a toy model, wich turns out to be
useful in the general description of how the algorithm acts in every
interation to get to an equilibrium price. The second one provides a
direct benchmark for the performance between our algorithm and a classical 
example in the literature, provided by Scarf \cite[Chapter 4]
{Kir98:elements_ge}. This section ends with a larger example of an
exchange economy (with symmetric agents), reflecting the computational
power of the augmented Walrasian approach. 

\begin{example}{\rm (symmetric agents).}\label{ex1.1}
To test the overall performance of the algorithm we start with a
basic example. Consider an economy of three goods and two agents, 
with utility functions within the CES family, i.e., 

\[ u_i(x)=\left(\sum_{j=1}^3 (a_{i,j})^{\frac{1}{b_i}}(x_j)^{\frac{b_i-1}{b_i}}\right)^{\frac{b_i}{b_i-1}}, \]

\noindent
with survival sets $X_i=[10^{-3},\infty)^2$, for each agent.
In this first example, the agents are symmetric, i.e., their utility functions' coefficients are equal, given by
 $a_{i,j}=\frac{1}{3}$, $i=1,2,\,j=1,2,3$ and $b_i=\frac{1}{2}$, $i=1,2$,
as well as their initial endowments $e_{i,j}=1$, $j=1,2,3\,,i=1,2$. 
It is easy to see that, by symmetry of the agents, the equilibrium 
price for this economy is given by $p^*=(\frac{1}{3},\frac{1}{3},\frac{1}{3})$,
and it is unique. Computationally, we initialize the algorithm at 
an arbitrary point of the simplex, in this case, $p^0=(0.12, 0.56, 0.32)$.
The trajectory of prices $\{p^\nu\}$ and excess supply evaluations $s(p^\nu)$ 
performed by our algorithm are depicted in Figure (\ref{ex1.1f}). The first 
graph describes the price evolution, where each good is represented by a line
(prices are scaled by a factor of 100). The second graph depicts the behaviour of 
the corresponding excess supply function, where each good is again represented by a line.
The adjustment process of the prices shows the \emph{Walrasian auctioneer's problem},
where in every iteration, the algorithm identifies the good of the excess supply
function with the least value, and performs an iteration adjusting its price 
for the next period. As the algorithm progresses, it converges to the equilibrium price.

% Example given by 20140430_2150
\begin{center}
\begin{figure}[ht]
%\scalebox{0.42}{\includegraphics{Numerical/AD/Results_20140430_2150/Prices_BW}}
\hfill
%\scalebox{0.42}{\includegraphics{Numerical/AD/Results_20140430_2150/Ex_Supply_BW}}
\caption{Homogeneus Agents (Example \ref{ex1.1})}\label{ex1.1f}
\end{figure}
\end{center}

\end{example}

%\begin{example}{\rm (heterogeneus agents).}\label{ex1.2}
%The second example corresponds to an heterogenous economy, considering
%the same type of utility functions (CES) with the following coefficients 
%and initial endowments

%\state Agent 1. $b=3.456$, $a=(0.419, 0.168, 0.413)$ and 
%$e=(0.926, 4.164, 4.910)$.

%\state Agent 2. $b=6.608$, $a=(0.429, 0.336, 0.235)$ and
%$e=(6.902, 0.533, 2.565)$.

%\noindent In this example, considering the same survival sets as the 
%previous one, initial price system is set to be the centroid of the simplex
%$p^0=(\frac{1}{3},\frac{1}{3},\frac{1}{3})$. The trajectory of prices 
%$p^k$ and excess supply evaluations $s(p^k)$ is depicted in Figure
%\ref{ex1.2f}.

%% Example given by 20140430_2129
%\begin{center}
%\begin{figure}[ht]
%\scalebox{0.42}{\includegraphics{Numerical/AD/Results_20140430_2129/Prices_BW}}
%\hfill
%\scalebox{0.42}{\includegraphics{Numerical/AD/Results_20140430_2129/Ex_Supply_BW}}
%\caption{Heterogeneus Agents (Ex.\ref{ex1.2})}\label{ex1.2f}
%\end{figure}
%\end{center}
%\end{example}

\begin{example}{\rm (exchange economy; Scarf example).}\label{ch1:scf}
Consider the example described in H. Scarf in {\rm \cite[Chapter 4]
{Kir98:elements_ge}}:
exchange economy involving five type of consumers and ten comodities. The initial 
endowment for each agent is given by Table \ref{ch1_scfe}.

\begin{table}
\label{ch1_scfe}
\caption{Initial endowments for Example \ref{ch1:scf}}
\begin{center}
\begin{tabular}{crrrrrrrrrr}
\hline
Consumer&\multicolumn{10}{c}{Initial endowments $e_{ij}$}\\
\hline
1&0.6&0.2&0.2&20.0&0.1&2.0&9.0&5.0&5.0&15.0\\
2&0.2&11.0&12.0&13.0&14.0&15.0&16.0&5.0&5.0&9.0\\
3&0.4&9.0&8.0&7.0&6.0&5.0&4.0&5.0&7.0&12.0\\
4&1.0&5.0&5.0&5.0&5.0&5.0&5.0&8.0&3.0&17.0\\
5&8.0&1.0&22.0&10.0&0.3&0.9&5.1&0.1&6.2&11.0\\
\hline
\end{tabular}
\end{center}
\end{table}

\noindent
The utility functions correspond to the CES-type, for which 
the parameters $a_{ij}$ and $b_i$ for each consumer are
described in Table \ref{ch1_scfut}.

\begin{table}[h]
\label{ch1_scfut}
\caption{Utility parameters for Example \ref{ch1:scf}}
\begin{center}
\begin{tabular}{crrrrrrrrrrr}
\hline
Consumer&\multicolumn{11}{c}{Utility parameters}\\
&\multicolumn{10}{c}{$a_{i,j}$}&$b_i$\\
\hline
1&1.0&1.0&3.0&0.1&0.1&1.2&2.0&1.0&1.0&0.07& 2.0\\
2&1.0&1.0&1.0&1.0&1.0&1.0&1.0&1.0&1.0&1.0& 1.3\\
3&9.9&0.1&5.0&0.2&6.0&0.2&8.0&1.0&1.0&0.2& 3.0\\
4&1.0&2.0&3.0&4.0&5.0&6.0&7.0&8.0&9.0&10.0& 0.2\\
5&1.0&13.0&11.0&9.0&4.0&0.9&8.0&1.0&2.0&10.0& 0.6\\
\hline
\end{tabular}
\end{center}
\end{table}

\noindent
The algorithm was set with the self-dual augmenting function,
and the centroid of the simplex as the initial point. Additionally,
the augmenting parameter is updated by $r^{\nu}=1.259$. The 
trajectory of prices $\{p^\nu\}$ and the corresponding sequence
of excess supply evaluations are depicted in Figure \ref{ex1.scf}.
In this example, the convergence to an approximate equilibrium point
for $\epsilon=10^{-1}$ is obtained within 37 iterations, taking a machine time
of 114 [min]; for $\epsilon=10^{-2}$, 53 iterations were required taking 179 [min]. 
The price is given by
$$p^*=\left(18.4,\,11.0,\,9.9,\,4.4,\,12.5,\,7.7,\,11.7,\,10.2,\,9.9,\,4.3\right)$$

\noindent
As in the previous example, the price sequence describes a trajectory that
can be associated with the \emph{Walrasian auctioneer's problem}. Similar
results are obtained with different starting points, as well as different
augmenting sequences.
%\begin{table}[h]
%\label{ch1:scf_b}
%\caption{Utility parameters $b$ for Example \ref{ch1:scf}}
%\begin{center}
%\begin{tabular}{cr}
%\hline
%Consumer&$b_i$\\
%\hline
%1&2.0\\
%2&1.3\\
%3&3.0\\
%4&0.2\\
%5&0.6\\
%\hline
%\end{tabular}
%\end{center}
%\end{table}

% Example given by Results_20141104_1951
\begin{center}
\begin{figure}[ht]
%\scalebox{0.42}{\includegraphics{Numerical/AD/Results_20150413_1557/Prices_BW}}
\hfill
%\scalebox{0.42}{\includegraphics{Numerical/AD/Results_20150413_1557/Ex_Supply_BW}}
\caption{Scarf's example (Example \ref{ch1:scf})}\label{ex1.scf}
\end{figure}
\end{center}

\end{example}

\begin{example}{\rm (large scale, symmetric agents economy).}\label{ex1.3}
In this example, we consider a larger economy, with  a total of 50 consumption goods and 10 agents with homogeneous
CES utility functions defined over survival sets given by $[10^{-3},\infty)^{50}$. The starting price 
is a random point in the simplex. As expected, the trajectory of the approximating
prices $\{p^\nu\}$ converges to the unique equilibrium price system, in which every 
good has the same price, i.e., $p_g=\frac{1}{50},\,g=1,\ldots,50$. The converge 
of the sequence of prices, $\{p^\nu\}$ and the corresponding sequence
of excess suppy functions $\{s(p^\nu)\}$ is illustrated in Figure \ref{ex1.3f}.
%Total time 11784.287185 [s]
% Example given by Results_20140501_0113
\begin{center}
\begin{figure}[ht]
%\scalebox{0.42}{\includegraphics{Numerical/AD/Results_20140501_0113/Prices_BW}}
\hfill
%\scalebox{0.42}{\includegraphics{Numerical/AD/Results_20140501_0113/Ex_Supply_BW}}
\caption{Large scale, symmetric agents (Example \ref{ex1.3})}\label{ex1.3f}
\end{figure}
\end{center}
\end{example}

From the examples previously described, a crucial observation can be made 
regarding the stability of the iterative process: the algorithm approaches
an approximating equilibrum with about half of the total iterations. This 
behaviour is robust in every simulation performed, and one find a reason
in the introduction of the augmenting function.

It's noteworthy that, in all cases, after a few iteration, the procedure finds
an approximate equilibrium which one should be able to exploit when dealing
with equilibirum problems in a stochastic environment.

% April 8, 2015 (Wednesday, 11:40 PDT) rjbw

\section{Dynamic deterministic equilibrium model}

As a stepping stone to the solution of stochastic Walras equilibrium
models, we are going to rely on solving, efficiently, deterministic
dynamic versions of the Walras equilibrium model. Our starting point is
a two-stage model that's formulated as follows: Given a price system
$p = (p^0, p^1)$ with $p^t$ the price vector in vigor at time $t$,
each agent $i\in \cI$ determines its optimal consumption plan
$\bar x = (\bar x_i^0, \bar x_i^1)$ as the
solution of the following utility maximization problem,
\begin{align*} 
\max_{x^0,y,x^1}\; &u_i^0(x^0) + u_i^1(x^1) \\
\sothat &\langle p^0, x^0 + T_i^0y\rangle \leq \langle p^0, e_i^0
\rangle, \\
 &\langle p^1, x^1 \rangle \leq \langle p^1, e_i^1 + T_i^1y\rangle, \\
 & \;\; x^0 \in X_i^0,\quad y\in Y_i, \quad x^1\in X_i^1,
\end{align*} 
where $u_i^t, e_i^t$ and $X_i^t$ are the utility functions, the
endowments and the survival sets for agent $i$ at time $t = 0,1$. 
As in \S 2, the utility functions are assumed to be usc and concave,
providing the convexity of the corresponding survivals sets. The
vector $y$ determines a set of activities selected by agent $i$ at time
$0$ that requires an input of goods $T_i^0y$ and produces a
deterministic output $T_i^1y$ at time $1$. The closed convex cone $Y_i
\subset \reals^m$ determines the set of potential activities that are at
the disposal of agent-$i$; in many instances one would simply have $Y_i
= \reals^m_+$ but not necessarily in general. One can think of the pair
of matrices $(T_i^0, T_i^1)$ as determining an input/output (home production)
process that could simply be savings including enhancements
or deterioration, or investment-activities, and so on. Of course, the agent
chooses $y$ so as to maximize its overall utility; note, $u_i^1$ could
include a discount factor that doesn't have to be made explicit here. \\

\noindent 
The excess supply function $s(p) = (s^0(p^0,p^1), s^1(p^0,p^1))$ is
given as usual as the difference between the total amount of goods
available in each time period and the total endowments adjusted by the
goods used or generated by the input/output process, i.e.,
\begin{align*} 
s^0(p) &= \sum\nolimits_{i\in \cI}[\,e_i^0 - (x_i^0(p) + T_i^0
y_i(p))\,], \\
s^1(p) &= \sum\nolimits_{i\in \cI}[\,(e_i^1+T_i^1 y_i(p)) - x_i^1(p)\,],
\end{align*} 
where $(x_i^0(p),y_i(p), x_i^1(p))$ is the optimal solution for agent
$i$ of its utility maximization problem. \\

\noindent 
The {\em Walrasian}, $W\cc\Delta^2\times\Delta^2\to \reals$ is defined by
\[W(p,q) = \langle q,s(p) \rangle = \langle (q^0,q^1),
		(s^0(p^0,p^1),s^1(p^0,p^1)) \rangle. \]
$\bar p = (\bar p^0, \bar p^1)$ is and equilibrium price system if
$s(\bar p) \geq 0$. As in the static (one-stage) model, it can be shown
that such a $\bar p$ is a maxinf-point of the Walrasian and its existence
is provided as $W$ is a Ky Fan function. One possible
approach in finding such a maxinf-point is based on the Augmented Walrasian
approach described in \S 3.

\begin{theorem}{\rm (dynamic deterministic maxinf-points).}\label{maxinf2stdet}
Consider the Walrasian function $W$ for the previous economy. Assuming local 
nonsatiation preferences, every maxinf-point $\bar{p}=(\bar{p}^0,\bar{p}^1)$ 
of $W$ is an equilibrium point, i.e., $s^0(\bar{p})\geq 0$ and $s^1(\bar{p})\geq 0$.
\end{theorem}

\state Proof.
Adapting the same pattern of proof as in Lemma \ref{lemm:maxinf}, for 
every price system $p=(p^0,p^1)$, $\langle p^0, s^0(p)\rangle=0$ and 
$\langle p^1, s^1(p)\rangle =0$. Then, if $\bar{p}$ is a maxinf-point of $W$, 
$W(\bar{p},\cdot)\geq 0$, and it follows that for vectors $q=(e^j,\bar{p}^1)$
defined for every unit vector $e^j$, $0\leq \langle q, s(\bar{p}) \rangle=\langle e^j,s^0(\bar{p})\rangle+\langle \bar{p},s^1(\bar{p})\rangle$ which implies $s^0(\bar{p})\geq 0$. 
Analogously, taking $q=(\bar{p}^0,e^j)$ it follows that $s^1(\bar{p})\geq 0$.
\eop

\begin{theorem}{\rm (convergence of $\epsilon$-maxinf points and equilibrium).}\label{conv_mi_2stdet}
Suppose that $p \mapsto s(p)$ is usc on $\Delta$. Consider the non-negative sequences 
$\lset r^\nu: \nu \in \nats \rset$ and $\lset \epsilon^\nu: \nu \in \nats \rset$ 
such that $r^\nu \upto \infty$, $\epsilon^\nu \downto \epsilon\geq 0$
\footnote{Note that the equilibrium case $\epsilon=0$ is included.}. 
Let $\lset W^{\nu}: \nu \in \nats \rset$ be a family of Augmented 
Walrasian functions associated which each augmenting parameter $r^\nu$. Let 
$p^{\nu} \in \epsilon^{\nu} \eqhyph \argmaxinf\, W^{\nu}$ 
and $\bar{p}$ be a cluster point of $\lset p^{\nu}:\nu \in \nats\rset$. 
Then $\bar{p} \in \epsilon \eqhyph \argmaxinf\, W$. In particular, when 
$\epsilon=0$, $\bar{p}$ is an equilibrium point.
\end{theorem}

\state Proof.
The tight lop-convergence of the augmented Walrasian
$\lset W^\nu : \nu \in \nats \rset$ follows the same arguments as those
cwin the proof of Theorem \ref{conv_mi_ad} and the 
conclusion follows from Theorem \ref{mot_3.2}.
\eop

\subsection{Dynamic model: a solution strategy}

For a fixed choice of activities $y_i$, the two-stage deterministic model
is essentially just an extension of a one-stage problem. With $y_i = \bar
y_i\in Y_i$, after dropping reference to agent-$i$, the problem reads:
\begin{align*} 
\max_{(x^0,x^1)}\; &u^0(x^0) + u^1(x^1) \\
\sothat &\langle p^0, x^0 \rangle \leq 
	\langle p^0, e^0 - T^0 \bar y\rangle, \\
&\langle p^1, x^1 \rangle \leq \langle p^1, e^1 + T^1 \bar y\rangle, \\
 & \quad x^0 \in X^0, \quad x^1\in X^1,
\end{align*} 
In fact, the problem is then separable, i.e., it can be solved 
by maximizing separately in the $x^0$ and $x^1$ variables:
\begin{align*} 
\nmax_{x^0\in X^0}\; &u^0(x^0)\; \sothat \langle p^0, x^0 \rangle \leq 
	\langle p^0, e^0 - T^0 \bar y\rangle, \\
\nmax_{x^1\in X^1}\; &u^1(x^1)\: \suchthat \langle p^1, x^1 \rangle \leq
	\langle p^1, e^1 + T^1 \bar y\rangle.
\end{align*} 
If these problems are of the Cobb-Douglas or CES-type, one can find
(closed-form) explicit solutions to these problems at negligible 
computational cost. \\

\noindent 
The agent's problem can now be seen as finding the best $y\in Y$ that will
maximize the overall rewards. With
\begin{align*}
 r(y) =& \sup_{x^0\in X^0} \lset u^0(x^0) \mset 
\langle p^0, x^0 \rangle \leq \langle p^0, e^0 - T^0 y\rangle \rset\\  
&+ \sup_{x^1\in X^1}  \lset u^1(x^1) \mset 
\langle p^1, x^1 \rangle \leq \langle p^1, e^1 + T^1 y\rangle \rset ,
\end{align*}
the agent's problem can be translated to:
\[\rwtxt{find}\; y^* \in \nargmax_{y\in Y} r(y). \]
We refer to this reduction as the \emph{transfer first approach} and
the algorithmic procedure to solve it (nonlinear convex
optimization problem) very much depends on the properties of $r$. In the
Cobb-Douglas or CES case, the function $r$ is twice
differentiable and one can find explicit expressions for the gradient
and the Hessian of $r$. When $Y= \reals_+^m$, the problem boils down to
maximizing a convex function on the non-negative orthant. Assuming
further that $r$ is differentiable, the optimality conditions read:
\[ \rwtxt{for}\; k = 1,\dots, m, \quad y_k^*\geq 0, \quad
    \frac{\partial}{\partial y_k}r(y^*) \leq 0, \quad 
    y_k^* \frac{\partial}{\partial y_k}r(y^*) = 0. \]
A number of specialized algorithmic procedures have been designed
for precisely this problem-type.

\subsection{The Cobb-Douglas case}

The utility function of agent-$i$ takes the form 
\[ u_i(x) = \nprod_{j=1}^n x_j^{\beta_{i,j}} \quad \rwtxt{with}
	\nsum_{j=1}\beta_{i,j} = 1, \;\; \beta_{i,j} \geq 0. \]
For $p\in\Delta$ and assuming that the survival set $X_i = \reals_+^n$,
agent-$i$ solution is
\[ \rwtxt{for}\; j = 1,\dots, n, \quad \bar x_{i,j}(p) =
\frac{\beta_{i,j}}{p_j} \nsum_{l=1}^n p_l e_{i,l}; \]
the endowment of agent-$i$: $e_i = (e_{i,1}, \dots, e_{i,n})$ and the
utility attached to this solution:
\[ u_i(\bar x_i) = \alpha_i(p) \Big(\nsum_{l=1}^n p_l e_{i,l}\Big) \quad
	\rwtxt{where}\; \alpha_i(p) = \nprod_{j=1}^N
	\Big(\frac{\beta_{i,j}}{p_j}\Big)^{\beta_{i,j}}. \]

\noindent For the dynamic model, once the activity levels $y \geq 0$
are {\em fixed}, the problem becomes {\em separable} (per-stage) and the
solution takes the same form {\em provided} that $y$ is chosen
so that $e_i^0 - T_i^0 y$ remains non-negative, otherwise agent-$i$
would enter the exchange market with a negative quantity of certain
goods. It's implicitly assumed that the technology matrices $T_i^0,
T_i^1$ are non-negative; negative entries in $T_i^0$ would imply
goods-production at time 0 and negative entries in $T_i^1$ would
imply negative outputs would be generated by certain technologies at
time 1. Hence, assuming that $T_i^0 y \leq e_i^0$, the solutions
(consumption vectors) that result from the choice of $y$ and
$p = (p^0,p^1)\in \Delta\times\Delta$ would be
\[ \rwtxt{for}\; j = 1,\dots, n, \quad \bar x_{i,j}^0(p^0) =
\frac{\beta_{i,j}^0}{p^0_j} \nsum_{l=1}^n p_l^0
	(e_{i,l}^0-\langle T_{i,l}^0, y \rangle); \]
where $T_{i,l}^0$ is the $l$th row of $T_i^0$,
\[ \rwtxt{for}\; j = 1,\dots, n, \quad \bar x_{i,j}^1(p^1) =
\frac{\beta_{i,j}^1}{p^1_j} \nsum_{l=1}^n p_l^1
	(e_{i,l}^1+\langle T_{i,l}^1, y \rangle); \]
and consequently,
\begin{align*} 
r_i(y) &= u_i^0(\bar x^0) + u_i^1(\bar x^1) \\
&= \alpha_i^0(p^0) \Big(\nsum_{l=1}^n p_l^0 (e_{i,l}^0 
		- \langle T_{i,l}^0, y \rangle)\Big) +
 \alpha_i^1(p^1) \Big(\nsum_{l=1}^n p_l^1 (e_{i,l}^1 +
 	\langle T_{i,l}^1, y \rangle)\Big) 
\end{align*} 
As detailed in \S 4.1, the optimization problem for agent-$i$ is reduced to
\[ \rwtxt{find}\; \bar y_i \;\; \rwtxt{that maximizes} r_i(y)
\;\suchthat
 T_i^0 y \leq e_i^0, \; y \in \reals_+^m. \]
This is a linear programming problem whose feasible region is bounded
and non-empty; $y = 0$ is always a feasible solution.

\subsection{The Constant Elastiticity of Substitution case.}

If the utility functions for agent $i$ take the following form
\begin{eqnarray*}
u_i^0(x^0) & = & \left( \sum_{j=1}^n ( a^0_{i,j} )^{\frac{1}{b_i^0}}(x^0_j)^{\frac{b_i^0-1}{b_i^0}} \right)^{\frac{b_i^0}{b_i^0-1}}.\\
u_i^1(x^1) & = & \left( \sum_{j=1}^n ( a^1_{i,j} )^{\frac{1}{b_i^1}}(x^1_j)^{\frac{b_i^1-1}{b_i^1}} \right)^{\frac{b_i^1}{b_i^1-1}}.
\end{eqnarray*}
Then, the KKT optimality conditions (\cite{VaAn},) are satisfied if, 
and only if, the budget constraint is active. On the other hand, each 
agent must satisfy the constraint for feasibility $T_i^0y\leq e_i^0$. 
Then, for a given a feasible $y\in Y_i$, we can find an explicit solution, 
given by

\[ \rwtxt{for}\; j = 1,\dots, n, \quad \bar x_{i,j}^0(p) =
\frac{ a^0_{i,j} }{ ( p^0_j )^{ b_i^0 } \sum_{ k=1 }^n ( p^0_k )^{ 1-b_i^0 } a_{i,k}^0 } \sum_{l=1}^n p^0_l ( e_{i,l}^0- \langle T_{i,l}^0 , y \rangle ); \]

\noindent where $T_{i,l}^0$ is the $l$th row of $T_i^0$,

\[ \rwtxt{for}\; j = 1,\dots, n, \quad \bar x_{i,j}^1(p^1) =
\frac{ a^1_{i,j} }{ ( p^1_j )^{ b_i^1 } \sum_{ k=1 }^n ( p^1_k )^{ 1-b_i^1 } a_{i,k}^1 } \sum_{l=1}^n p^1_k ( e_{i,l}^1 + \langle T_{i,l}^1 , y\rangle ); \]

\noindent Defining for agent-$i$
\[ \rwtxt{for}\; t = 1,2, \quad\theta_i^t(p) = \Big( \sum_{j=1}^n ( a^t_{i,j} )^{ \frac{ 1 }{ b_i^t } } \Big( \frac{ a_{i,j}^t }{ ( p^t_j )^{ b_i^t }} \frac{  1 }{ \sum_{k=1}^n ( p^t_k )^{ 1-b_i^t }a_{i,k}^t } \Big)^{ \frac{ b_i^t-1 }{ b_i^t } } \Big)^{ \frac{ b_i^t }{ b_i^t-1 } }. \]

\noindent consequently
\begin{align*} 
r_i(y) &= u_i^0(\bar x^0) + u_i^1(\bar x^1) \\
&= \theta_i^0(p^0) \Big(\nsum_{l=1}^n p_l^0 (e_{i,l}^0 
		- \langle T_{i,l}^0, y \rangle)\Big) +
 \theta_i^1(p^1) \Big(\nsum_{l=1}^n p_l^1 (e_{i,l}^1 +
 	\langle T_{i,l}^1, y \rangle)\Big) 
\end{align*} 

\noindent This is a linear function of $y$. Thus, if $Y_i=\reals_+^m$, the 
problem for each agent is given by
\[ \rwtxt{find}\; \bar y_i \;\; \rwtxt{that maximizes} r_i(y)
\;\suchthat
 T_i^0 y \leq e_i^0, \; y \in \reals_+^m. \]

\section{Stochastic Equilibrium.}

\subsection{The agent's problem.}

In an uncertain (stochastic) environment the agent's problem 
two-stage problem can be formulated as follows:
\begin{align*} 
\max_{x^0,y,x_{\cdot}^1}\; &u_i^0(x^0) + E_i\{u_i^1(\bfxi,x_{\bfxi}^1)\} \\
\sothat &\langle p^0, x^0 + T_i^0y\rangle \leq \langle p^0, e_i^0 \rangle, \\
 &\langle p_{\xi}^1, x_{\xi}^1 \rangle \leq \langle p_{\xi}^1,
   e_{i,\xi}^1 + T_{i,\xi}^1 y\rangle, \quad\forall\,\xi\in\Xi\\
 & y\in \reals_+^m, \;\; x^0 \in X_i^0, \;\; x_{\xi}^1\in X_{i.\xi}^1,
   \quad\forall\,\xi\in\Xi,
\end{align*}
where, the utility functions are usc, concave and the survival 
sets are convex and unbounded. Additionally, the set $\Xi$ 
consists of a finite number of possible states (scenarios)
and $E_i\{\cdot\}$ indicates that agent-$i$ is calculating the
expectation with respect to agent-$i$ {\em beliefs}, i.e., to each possible
state $\xi\in\Xi$, agent-$i$ assigns a probability $\pi_{i,\xi} \geq 0$
such that $\sum_{\xi\in\Xi} \pi_{i,\xi} = 1$.  It's possible, 
although unlikely, that all agents have the same information 
about the future in which case these probabilities wouldn't 
depend on $i$. As before, the agents set up their trades in full 
knowledge of the suggested price system, eventually an equilibrium
price system,
\[ p = (p^0, (p_\xi^1)_{\xi\in\Xi} ); \]
in particular, $p_\xi^1$ is known for {\em every} contingency
$\xi\in\Xi$. Note that the goods required $T_i^0 y$ to carry
out activities at level $y$ are still well determined, the output at time
1 is now stochastic, namely $T_{i,\xi}^1y$. This reflects a more realistic
view of the output process. Even in the simple case of savings via buying
certificates of deposit, bonds or stocks, their value at time 1 can't be
known with certainty. This is even more so, if the activities are
decisions involving manufacturing, the marketing or distribution
of goods (perishable or not)  and so on.\\

\noindent The agents' problems are thus (two-stage)
{\em stochastic programs with recourse} \cite{BL2011} with stochastic
entries in the right-hand side $\langle p_{\xi}^1,
e_{i,\xi}^1\rangle$, the so-called technology matrix $
T_{i,\xi}^\top$ and the recourse matrix $p_{\xi}^1$; the
recourse decisions are $x_{i,\xi}^1$. Under the `usual' conditions that guarantee
the existence of an equilibirum price system, recalled in \cite{JfrW02:wlrs, MQ2002},
these stochastic programs are necessarily feasible; note however that straightforward
feasibility of these stochastic programs doesn't really require such stringent
conditions, for example, one could rely on an adaptation of the {\em ample
survivavibility assumption} introduced in \cite{JfRW14:gei}.
From the stochastic programming viewpoint these conditions can be viewed as
sufficient conditions to guarantee the relatively complete recourse property. 
For our problem, this can be stated as follows: for every agent,
(dropping the dependence on $i$), for all $\xi\in\Xi$, there exists 
$(\tilde{x}^0,\tilde{y},\tilde{x}^1_\xi)\in X^0\times \reals^m_+\times X^1_\xi$
such that 
\begin{eqnarray*}
e^0_l-\tilde{x}^0_l-(T^0\tilde{y})_l & \geq & 0,\quad l=1,\ldots,L\\
e^1_{l,\xi}-\tilde{x}^1_{l,\xi}+(T^1_\xi\tilde{y})_l & \geq & 0,\quad l=1,\ldots,L,\,\forall \xi\in\Xi
\end{eqnarray*}

\noindent From the economic perspective, this assumption
is weaker than the usual survability conditions, and can 
be interpreted as every agent being able to survive or
participate in the economy, independently of the market prices.

\subsection{Solving the agent's (stochastic) problem}

There are many alternatives methods to solve stochastic programs
with recourse,
but in this setup the use of the Progressive Hedging algorithm
\cite{RckW91:scn, Wets89:aggreg} seems to
have many advantages, in particular because solutions of the
individual scenario subproblems are so readily available, cf.\ \S 4 with
the Cobb-Douglas case. \\

\noindent The approach is based on relaxing, at the outset, the {\em non-anticipativity}
constraint, namely that $x^0$ and $y$ aren't allowed to depend on $\xi$, and
then, {\em progressively} enforcing this requirement.  For now let's
just limit ourselves to a description of the steps of the algorithm as
it applies to the stochastic version of the (two-stage) agent's problem
in the Cobb-Douglas case, generated by the \emph{transfer first} approach
described in \S 4. \\

\noindent
{\bf Step 0}. Set $\nu=0$. Pick $\rho > 0$, $\bar y^0 = 0$,
$w_i^\nu: \Xi\to \reals^m$ 
such that $E_i\{w_i^\nu(\bfxi)\} = 0$. \\

\noindent 
{\bf Step 1}. For all $\xi\in\Xi$, let
\[ y_i^{\nu+1}(\xi) \in \nargmax_y \lset r_i^\nu(\xi,y) -
   \langle w_i^\nu,y\rangle - \frac{\rho}{2}|y - \bar y^\nu|^2
   \mset T_i^0 y \leq e_i^0, \;\; y\in\reals_+^m \rset,\]
where 
\begin{align*}
r_i^\nu(\xi,y) 
= &\alpha_i^0(p^0) \Big(\nsum_{l=1}^n p_l^0 (e_{i,l}^0 
   - \langle T_{i,l}^0, y \rangle)\Big)\\
 &+ \alpha_i^1(\xi,p^1(\xi)) \Big(\nsum_{l=1}^n p_l^1(\xi) (e_{i,l}^1(\xi) +
 	\langle T_{i,l}^1(\xi), y \rangle)\Big)
\end{align*}
and 
\[\alpha_i^0(p^1) = \nprod_{j=1}^n
\Big(\frac{\beta_{i,j}^0}{p_j^0}\Big)^{\beta_{i,j}^0},
\qquad 	\alpha_i^1(\xi,p^1(\xi)) = \nprod_{j=1}^n
\Big(\frac{\beta_{i,j}^1}{p_j^1(\xi)}\Big)^{\beta_{i,j}^1}. \]

\noindent
{\bf Step 2}. If $\xi \mapsto y_i^\nu(\xi)$ is a constant function, stop.
$y^\nu(\xi)$, for any $\xi$, of course, determines the optimal activity
levels and the corresponding vector and function 
$[x_i^0, (x_i^1(\xi), \xi\in\Xi)]$ determine the optimal consumption
plans.
Otherwise, set $\bar y_i^{\nu+1} = E\{y_i^{\nu+1}(\bfxi)\}$,
\[w_i^{\nu+1}(\xi) = w_i^\nu + \rho(y_i^{\nu+1}(\xi) - \bar
y_i^{\nu+1}), \]
and return to {\bf Step 1} with $\nu = \nu+1$. \\

\noindent 
Note that the optimization problem in {\bf Step 1} is a quadratic
program of a very simple nature since it's completely separable. After
carrying out some elementary calculations, it can be written in the form:
\[ \max \lset \nsum_{j=1}^m (\bar c_j(\xi) y_j - \frac{\rho}{2}y_j^2) \mset 
       T_i^0 y \leq e_i^0, \;\; y\in\reals_+^m \rset.\]
One could rely on general quadratic procedures to solve this particular
problem, but a much more efficient procedure could be designed to
deal with a problem of this particular type.\\

%%%%%%%%%%%%%%%%%%%%%%%%%%%%%%%%%%%%%%%%%%% jd
% Should we put the description of the Augmented Walrasian here or 
% before the solution strategy of the Agent's problem?

\noindent One final remark about this model is that it can be easily extended
to the CES utility functions case, using the same transfer first approach of
maximizing $r$ function. In this situation, is easy to see that the only 
difference is the sustitution of the linear coefficientes $\alpha$ by
the ones given by the CES parameters, $\theta$. Furthermore, one can solve
the general agent problem, relaxing the dependence of $x^0$ and $y$ on $\xi$,
and apply the enforcing procedure to progressively converge to a 
\emph{deterministic} solution.

\subsection{Augmented Walrasian and approximating scheme.}
In this section, we set the foundations of the augmentation techniques applied to
the Dynamic Stochastic Equilibrium Model. A description of equilibrium points as
maxinf points of the corresponding Walrasian, as well as the approximation scheme 
based in tight lopsided convergence of augmented Walrasian is provided. Finally, 
a general description of the computational implementation of the algorithm
and numerical examples are analized.\\

\noindent 
As in the previous section, consider the {\em stochastic equilibrium model}, 
where given a price system $p=\left(p^0, (p_\xi^1)_{\xi\in\Xi}\right)$, each agent $i$ solves

\begin{align*} 
\max_{x^0,y,x_{\cdot}^1}\; &u_i^0(x^0) + E_i\{u_i^1(\bfxi,x_{\bfxi}^1)\} \\
\sothat &\langle p^0, x^0 + T_i^0y\rangle \leq \langle p^0, e_i^0 \rangle, \\
 &\langle p_{\xi}^1, x_{\xi}^1 \rangle \leq \langle p_{\xi}^1,
   e_{i,\xi}^1 + T_{i,\xi}^1 y\rangle, \quad\forall\,\xi\in\Xi\\
 & y\in \reals_+^m, \;\; x^0 \in X_i^0, \;\; x_{\xi}^1\in X_{i.\xi}^1,
   \quad\forall\,\xi\in\Xi,
\end{align*}

\noindent
which defines the individual demand function $ x_{i}(p) = \left(x^0_{i}(p), (x_{i,\xi}^1(p))_{\xi\in\Xi}\right)$ and the individual transfer vector $y_i(p)$. 
Additionally, the excess supply function for this economy $s(p) = \left(s^0(p),(s_\xi^1(p))_{\xi\in\Xi}\right)$ is defined as

\begin{eqnarray*}
s^0(p)&=&\sum_{i\in\mathcal{I}} e_i^0-x_i^0(p)-T^0_iy_i(p)\\
s^1_\xi(p)&=&\sum_{i\in\mathcal{I}} e_{i,\xi}^1-x_{i,\xi}^1(p)+T^1_{i,\xi}y_i(p), \quad\forall\,\xi\in\Xi\\
\end{eqnarray*}

\noindent
The Walrasian for this model is the function $W:(\Delta\times\Delta^{|\Xi|})\times(\Delta\times\Delta^{|\Xi|})\to\reals$
defined by

\[ W(p,q)  = \langle q^0, s^0(p)\rangle+\sum_{\xi\in\Xi} \langle q_\xi^1, s_\xi^1(p)\rangle 
\]

\noindent
A price system $\bar{p}=\left(\bar{p}^0, (\bar{p}_\xi^1)_{\xi\in\Xi}\right)$ is an 
equilibrium price if $s(\bar p)\geq 0$, i.e., $s^0(\bar{p})\geq 0$ and $s_\xi^1({p})\geq 0$, for
every possible state $\xi\in\Xi$. Then, a equilibrium point for the
dynamic stochastic model can be described as a maxinf point of the Walrasian. Again, 
the existence is granted by noting that $W$ turns out to be a Ky Fan function.

% insert comment about existence of equilibrium points.

\begin{theorem}{\rm (stochastic equilibrium prices and maxinf-points).}\label{2ststoeq}
Consider the Walrasian function $W$ for the previous economy. Then, under local 
nonstatiation of preferences, every maxinf-point $\bar{p}=(\bar{p}^0,(\bar{p}_\xi^1)_{\xi\in\Xi})$ 
of $W$ is an equilibrium point, i.e., $s^0(\bar{p})\geq 0$ and $s_\xi^1({p})\geq 0$, for
every possible state $\xi\in\Xi$.
\end{theorem}

\state Proof.
Considering that for every price system $p=(p^0,(p^1_\xi)_{\xi\in\Xi})$, under local
nonsatiation preferences, the excess supply satisfies the Walras' law for the first stage
and for every possible state of the second stage, i.e.,$\langle p^0, s^0(p)\rangle =0$
and for every $\xi$, $\langle p^1_\xi, s^1_xi(p)\rangle =0$. Thus, for $\bar{p}$ a maxinf 
point of $W$, $W(\bar{p},\cdot)\geq 0$. Considering $q=(e^j,(\bar{p}^1_\xi)_{\xi\in\Xi})$, 
$0\leq \langle q, s(\bar{p})\rangle=\langle e^j,s^0(\bar{p})\rangle + \sum_{\xi\in\Xi}\langle \bar{p}^1_\xi), s^1_\xi(\bar{p})\rangle$, which implies that $(s^0(\bar{p}))_j\geq 0$ for all $j$. For the second stage, given an scenario $\xi_0\in \Xi$, it suffices to take $q=(\bar{p}^0,(\bar{p}^1_1,\ldots,\bar{p}^1_{\xi_0-1}, e^j, \bar{p}^1_{\xi_0+1},\ldots,\bar{p}^1_{\Xi})$ and conclude by the same argument that $(s^1_{\xi_0}(\bar{p}))_j\geq 0$, for all $j$ and all $\xi_0$.
\eop

\noindent 
For the problem of finding equilibrium points, we will follow the 
approximating technique  described in \S 3, where for a given 
Walrasian function $W$ for the stochastic economy, 
we consider the  augmenting function $\sigma$, and an increasing 
sequence of positive scalars $r^\nu\upto\infty$, for which we defined 
the family of augmented Walrasian  bifunctions $W^\nu$ as follows 

\[ W^\nu(p,q) = \inf_{z\in\Delta\times\Delta^{\Xi}}\lset W(p,z)+r^\nu\ast \sigma^*(q-z)\rset, \]

\noindent
and the algorithmic procedure relies in the idea of finding 
approximating maxinf-points of this augmented Walrasian bifunctions 
for $\nu$ large enough. Finally, the following convergence result will
guarantee the approximation to an equilibrium point for the initial
economy. 

\begin{theorem}{\rm (convergence of dynamic stochastic $\epsilon$-maxinf points).}\label{conv2ststo}
Suppose that $p \mapsto s(p)$ is usc on $\Delta$. Consider the non-negative sequences 
$\lset r^\nu: \nu \in \nats \rset$ and $\lset \epsilon^\nu: \nu \in \nats \rset$ 
such that $r^\nu \upto \infty$, $\epsilon^\nu \downto \epsilon$, for $\epsilon\geq0$. 
Let $\lset W^{\nu}: \nu \in \nats \rset$ be a family of Augmented 
Walrasian functions associated wich each augmenting parameter $r^\nu$. Let 
$p^{\nu} \in \epsilon^{\nu} \eqhyph \argmaxinf\, W^{\nu}$ 
and $\bar{p}$ be a cluster point of $\lset p^{\nu}:\nu \in \nats\rset$. 
Then $\bar{p} \in \epsilon \eqhyph \argmaxinf\, W$. In particular, for 
$\epsilon=0$, $\bar{p}$ is an equilibrium point.
\end{theorem}

\state Proof.
The proof follows from the application of the Theorem \ref{mot_3.2}, as it was used in
the convergences results of sections \S 3 and \S 4 (Theorem \ref{conv_mi_ad}, Theorem \ref{conv_mi_2stdet}). 
Finally, the tight lopsided convergence of the sequence $\lset W^\nu: \nu\in\nats\rset$
follows from the same argument.
%: for $q\in\Delta\times\Delta^\Xi$, $\lset p^\nu:\nu\in\nats\to p\in\Delta\times\Delta^\Xi$ and defining $q^\nu:=1,\,\nu\in\nats$, then
%\[
%W^\nu(p,q)  =  \inf_{z\in\Delta\times\Delta^\Xi} \lset W(p,q)+r^\nu\ast \sigma^*(q^\nu-z)\rset \leq W(p^\nu,q^nu).
%\]
%Under usc of the function $p\mapsto s(p)$, $\limsup W^\nu(p^\nu,q^\nu)\leq \limsup W^\nu(p^\nu,q)\leq  W(p,q)$.  On the other hand, the function $F^\nu(q,z):=W(p,q-z)+r^\nu\ast \sigma^*(q-z)$ is
%level bounded in $z$ locally uniform in $q$. Therefore $W^\nu(p,q)$ is the inf-projection
%of $F^\nu$ in the $z$-variable and, in virtue of theorem \cite[Thm.1.17]{VaAn}, lsc. Finally,
%for any $q_0\in\Delta\times\Delta^\Xi$, $W^\nu(p,q^0)\to W(p,q^0$, and by a diagonal argument,
%$\liminf W^\nu(p^\nu,q^\nu)\geq W(p,q)$.
\eop

\subsection{Numerical implementation and examples.}
\noindent
Computationally, we proceed with a primal-dual iteration scheme as 
described in \S 3. Especial features for this type of economy are
considered. In terms of the agent's problem, we can adopt a strategy
solving the problem directly or solving it through the maximization
of the overall reward function $r$.\\

\noindent On the other hand, the agent's problem is a 
stochastic program with relatively complete recourse, for which 
Progressive Hedging algorithm is implemented. Exploiting the 
structure of the agent's problem given by the separability 
in terms of the different scenarios in the second stage, combined with the 
progressive hedging approach, we provide two strategies, one 
sequential and another one parallel. The efficiency of these strategies
will  be discussed later and will basically depend on the size of 
the economy considered as the total amount of goods available.\\

\noindent Finally, the global strategy of solution adopted can be 
summarized in the following scheme:

\noindent
{\bf Step 0}. Set $\nu=0$. Pick an initial price $p^{(0)}$ (for example, the centroid of the $\reals^G$-simplex for each $\xi\in\Xi$), and an augmenting parameter $r^0>0$. Define an strategy for the agent's problem, directly maximizing the utility function $u(x^0,y,(x^1_\xi)_{\xi\in\Xi})$, or indirectly maximizing the overall reward
    function $r(y)$. Additionally, defined the procedure for the Progressive Hedging 
    algorithm implementation, sequential or parallel.\\

    \noindent 
    {\bf Step 1}. For all $i\in\mathcal{I}$, compute $x_i^\nu(p^\nu)$ applying Progresive Hedging algorithm to the agent's problem with the proper choice of strategies. With this, compute $s^\nu(p^\nu)$, and solve the Phase I iteration for the primal-dual
    scheme:
    \[ q^{\nu+1} \in \nargmax_q \lset W_{r^\nu}(p^\nu,q) 
     \mset q\in\Delta\times\Delta^{\Xi} \rset,\]
    which is a linear problem.\\

    \noindent 
    {\bf Step 2}. Solve the Phase II, given by
    \[ p^{\nu+1} \in \nargmax_p \lset W_{r^\nu}(p,q^{\nu+1}) 
     \mset p\in\Delta\times\Delta^{\Xi} \rset.\]

    \noindent
    Finally, check the optimality condition: if $\min s(p^{\nu+1})\geq - \varepsilon$, 
    stop. Otherwise, set $r^{\nu+1}>r^\nu$ and return to {\bf Step 1} 
    with $\nu = \nu+1$. \\

    \subsection{Numerical experimentation.}

\begin{example}{\rm (main example).}\label{ex5.1}
The main example testing the numerical implementation of the 
augmented Walrasian algorithm is described for an economy consisting of 
seven goods: skilled job, unskilled job, leisure, consumption, risk free bond, 
and two stocks. We considered an economy with five agents, with utility
functions of CES type, and nine posible scenarios in the second stage. On
the other hand, the transformation matrices are the same for every 
agent at the first stage given by $T^0={\bf I}$ and for the second stage are given by $T_{i,\xi}^1= diag(d_{i,\xi})$ for each agent $i=1,\ldots\cI$, with
\begin{eqnarray*}
d_{1,\xi}&=&(0,0,1+3r/4,0.7,1+r,R^1_\xi,R^2_\xi),\\
d_{2,\xi}&=&(0,0,1+r/2,0.8,1+r,R^1_\xi,R^2_\xi),\\ 
d_{3,\xi}&=&(0,0,0,0.7,1+r,R^1_\xi,R^2_\xi),\\
d_{4,\xi}&=&(0,0,1+r/2,0.9,1+r,R^1_\xi,R^2_\xi),\\
d_{5,\xi}&=&(0,0,1+r/2,0.7,1+r,R^1_\xi,R^2_\xi).
\end{eqnarray*}

\noindent
where $r=3.25\%$ and $R^1_\xi$,$R^2_\xi$ are given by the following table
\begin{center}
\begin{tabular}{cc|ccc}
$\xi$&&$R^2_{(+)}$&$R^1_{(=)}$&$R^1_{(-)}$\\
&&1.10&1.00&0.95\\
\hline
$R^1_{(+)}$&1.20&1&2&3\\
$R^1_{(=)}$&1.00&4&5&6\\
$R^1_{(-)}$&0.85&7&8&9
\end{tabular}
\end{center}

\noindent
Agents' utility functions are CES type, with parameters can be found online \footnote{\url{http://www.math.ucdavis.edu/~jderide/AugWal/AugWal.html}},
as well as their initial endowments and survival sets. Additionally, 
we consider that every agent has the same beliefs over the scenarios
on the second stage, given by $\pi_{i,\xi}=\frac{1}{9},\,i\in\cI,\, \xi\in\Xi$.

\noindent
The algorithm is initialised with $p^{(0)}$ as the centroid of 
$\Delta\times\Delta^\Xi$, the augmenting function is 
$\sigma=\frac{1}{2}|\cdot|^2$, and the augmenting sequence of 
parameters $r^\nu$ is given by $r^\nu=1.259^\nu$.
The trajectory of the prices $\{p^\nu\}$ for every iteration 
and the corresponding excess supply function $\{s(p^\nu)\}$ are 
described in figure \ref{ex5.1f}. The algorithm was set for direct
solution for the agent's problem and for Progressive Hedging, a 
sequential approach was considered. It finished after 62 iterations,
with a total machine time of 28 [hrs]. 

% Example given by 20140430_2150
\begin{center}
\begin{figure}[ht]
%\scalebox{0.42}{\includegraphics{Numerical/2StSto/Results_20140725_0043/p0_BW}}
\hfill
%\scalebox{0.42}{\includegraphics{Numerical/2StSto/Results_20140725_0043/s0_BW}}
\caption{Main example (Example \ref{ex5.1}), $\{p^\nu\}$ and $\{s(p^\nu)\}$ }\label{ex5.1f}
\end{figure}
\end{center}
\end{example}

%\begin{example}{\rm (main example, parallel implementation).}\label{ex5.2}
%\noindent
%The parameters for this economy are the same considered in the previous
%example \ref{ex5.1}, but here we present the implementation considering
%a parallel strategy for the Progressive Hedging algorithm. The overall time 
%was 44 [hrs] with a total of 78 iterations. The trajectories for prices and 
%corresponding excess supply are described in figure \ref{ex5.2f}. 
%
%\begin{center}
%\begin{figure}[ht]
%\scalebox{0.42}{\includegraphics{Numerical/2StSto/Results_20140723_1544/p0_BW}}
%\hfill
%\scalebox{0.42}{\includegraphics{Numerical/2StSto/Results_20140723_1544/s0_BW}}
%\caption{Main example (Example \ref{ex5.1}), $\{p^\nu\}$ and $\{s(p^\nu)\}$ }\label{ex5.2f}
%\end{figure}
%\end{center}
%
%\noindent 
%For this economy, the total time of our parallel implementation is larger than
%using our sequential implementation, which can be explained by the size
%of the agent's problem not being large enough to justify the use of 
%a parallel strategy. 
%%A comparison between the iteration times for the  two cases is provided in figure \ref{ex5.2it}

%\end{example}

\section{Conclusions}

We introduced a new optimization methodology that allows the computation of 
equilibrium demand and prices for different economies. This 
new approach combines several elements of variational analysis, 
such as the notion of lopsided convergence and augmented 
Lagrangian  technique for non-concave optimization problems. 

Following \cite{JfrW02:wlrs}, we characterize equilibrium 
prices as maxinf points for the so-called Walrasian bifunction 
for an exchange economy. The novelty of our approach relies 
in the approximation of the Walrasian by \emph{augmented Walrasian}. 
Then, the computation of equilibrium points follows from the 
convergence of the sequence of maxinf points for the approximated problems.

We use this methodology to solve, as a prelude, the classical Arrow-Debreu 
general equilibrium model and, then, two periods exchange economies with 
uncertainty. For both models we got convergence in every numerical example, 
including a large scale problem in the stochastic case. A robust performance 
of the algorithm is always obtained, and it can be interpreted as a direct 
result of the augmentation procedure. One can appreciate stability of the 
iterations: by about half of the total iterations required to get a high 
tolerance-level solution. Furthermore, different numerical scenarios were 
tested, varying the augmenting function $\sigma$ and the augmenting 
parameter $r$. The results observed in these variations were not considered 
significantly different. The most efficient variant relied on the self-dual 
augmenting function with exponential growth in the augmenting parameter. 
Finally, for the stochastic problem, we tested an implementation of the 
algorithm based on a parallel computation for the agent problem. 

The usage of the augmented Walrasian approximation for the computation 
of equilibrium points can be extended for more sophisticated economic 
models, as the one presented in \cite{JfRW14:gei}, where financial markets, 
collateral, and retention goods are considered. Additionally, considering  
the structure of the problems, computational strategies that consider an 
efficient use of a parallel algorithm should improve the overall time performance. 

\subsection*{Acknowledgement}

This material is based upon work by Julio Deride and Roger Wets supported 
in part by the U.S. Army Research Laboratory and the U.S. Army Research Office 
under grant numbers W911NF-10-1-0246 and W911NF-12-1-0273. We thank David Woodruff,
University of California Davis, for his assistance with Pyomo, the optimization modeling
language used in the computational implementation of our algorithm.

%%%%%%%%%%%%%%%%%%%%%% jd
%\input{finMrkt}
%%%%%%%%%%%%%%%%%%%%%% jd
\bibliography{07-awlrs}

\end{document}